\documentclass{amsart}

\usepackage{amsmath,amsthm}
\usepackage{amssymb, enumitem, mathrsfs, verbatim, bm, esint}
\usepackage{amstext}
\usepackage[top=3.5cm,bottom=3.5cm,left=3.5cm,right=3.5cm]{geometry}

\usepackage{hyperref}

\newtheorem{theorem}{\bf Theorem}[section]

\newcommand{\N}{\mathbb N}

\newcommand{\R}{{\mathbb R}}
\newcommand{\ep}{\varepsilon}
\newcommand{\dual}[2]{\langle #1,#2\rangle}

\numberwithin{equation}{section}
\newcommand{\diver}{\operatorname{div}}
\DeclareMathOperator*{\essinf}{ess\,inf}
\DeclareMathOperator*{\sign}{sign}
\newcommand{\dert}{\partial_t}
\newcommand{\cni}{c_\nu}
\newcommand{\thetah}{\hat\thet}
\newcommand{\thetat}{\tilde\thet}

\newcommand{\vc}[1]{{\bf #1}}
\newcommand{\vv}{\vc{v}}

\newcommand{\vw}{\vc{w}}
\newcommand{\thet}{\vartheta}

\renewcommand{\S}{\mathbf{S}}
\newcommand{\D}{\mathbf{D}}

\newcommand{\dt}{\,{\rm d} t }
\newcommand{\dx}{\,{\rm d} {x}}
\newcommand{\dxdt}{\dx  \dt}
\newcommand{\dtau}{\,{\rm d} \tau}

\newcommand{\intO}[1]{\int_{\Omega} #1 \ \dx}
\newcommand{\intTO}[1]{\int_0^T \int_{\Omega} #1 \ \dxdt}

\newcommand{\nm}[1]{\|#1\|}

\begin{document}

\title{On solutions for a generalized Navier-Stokes-Fourier system fulfilling entropy equality}
\author[A.~Abbatiello]{Anna Abbatiello}
\address{Anna Abbatiello.
Sapienza University of Rome, Department of Mathematics ``G. Castelnuovo", Piazzale Aldo Moro 5, 00185 Rome, Italy.}   \email{\tt anna.abbatiello@uniroma1.it}

\author[M.~Bul\'{i}\v{c}ek]{Miroslav Bul\'i\v{c}ek}
\address{Miroslav Bul\'i\v{c}ek. Mathematical Institute, Faculty of Mathematics and Physics, Charles University, Sokolovsk\'{a} 83, 18675 Prague, Czech Republic.}
\email{\tt mbul8060@karlin.mff.cuni.cz}
 
\author[P.~Kaplick\'{y}]{Petr Kaplick\'{y}}
\address{Petr Kaplick\'{y}. Charles University, Faculty of Mathematics and Physics, Department of Mathematical Analysis, Sokolovsk\'{a} 83, 18675 Prague, Czech Republic.} 
 \email{\tt kaplicky@karlin.mff.cuni.cz}

%%%% Subject entries to be placed here %%%%

%%%% Keyword entries to be placed here %%%%

\begin{abstract}
We consider a flow of non-Newtonian heat conducting incompressible fluid in a bounded domain subjected to the homogeneous Dirichlet boundary condition for the velocity field and the spatially inhomogeneous Dirichlet boundary condition for the temperature. The ultimate goal is to show that the fluid converges to equilibrium as time tends to infinity. However, to justify such result, one needs to deal with very special inequalities and very special test functions, which are typically not admissible on the level of weak solutions. In this paper, we show how one can overcome such difficulties. In particular, we show the existence of a solution fulfilling the entropy equality, which seems to be optimal class of solutions in which one should study the stability result.
\end{abstract}

\maketitle

%\section{Formulation of the problem}
\section{Formulation of the problem}\label{S1}

We study the generalized Navier--Stokes--Fourier system
\begin{subequations}\label{bl}
\begin{align}
\label{i1}
\dert \vv + \diver(\vv \otimes \vv ) - \diver \S + \nabla \pi &= 0 \\
\label{i2}
\diver \vv &= 0\\
\label{i3}
\dert(\cni\thet)+\diver(\cni\thet\vv) + \diver \vc{q} &= \S:\D{\vv}
\end{align}
\end{subequations}
in $Q:=(0,T) \times \Omega\subset(0, +\infty)\times \R^3$ with $\Omega$ a bounded domain. The system~\eqref{i1}--\eqref{i3} is completed by the following boundary conditions
\begin{equation}
\vv=0,
\ \
\thet = \thet_{b} \qquad \mbox{ on } (0, T)\times\partial\Omega,
\end{equation}
where the function $\thet_{b}=\thet_b(x)$ is a nontrivial function of the position, and the following initial conditions
\begin{equation}\label{condizioni-iniziali}
\vv = \vv_0, \ \
\thet = \thet_0 \qquad \mbox{ in } \{0\}\times\Omega.
\end{equation}

Here $\vv: Q\to \R^3$ denotes the velocity field,  $\D\vv:=  (\nabla\vv + (\nabla\vv)^{t})/2$ is the symmetric part of the velocity gradient $\nabla \vv$, $\pi:Q\to \R$ is the pressure, $\thet:Q\to \R$ is the temperature; $\S: Q\to \R^{3\times 3}_{\rm sym}$ denotes the viscous part of the Cauchy stress tensor and $\vc{q}:Q \to \R^3$ is the heat flux.

Concerning the material parameters, the constant $\cni>0$ in \eqref{i3} denotes the heat capacity and, for simplicity and without losing the generality, we set $\cni\equiv 1$ in what follows. The heat flux $\vc{q}$ is represented by the Fourier law
\begin{equation}\label{Fourier}
\vc{q} = - \kappa(\thet) \nabla \thet
\end{equation}
with the heat conductivity $\kappa: \R \to (0, +\infty)$ being a continuous function of the temperature satisfying, for all $\thet \in (0, +\infty)$ and for some $0<\underline{\kappa}, \overline{\kappa} <+\infty,$
\begin{align}\label{k}
0<\underline{\kappa}\leq \kappa(\thet) \leq \overline{\kappa}<+\infty.
\end{align}
We assume that $\S=\S^*(\thet, \D\vv)$, where $\S^*:(0,\infty)\times \R^{3\times 3}_{\rm sym} \to \R^{3\times 3}_{\rm sym}$ is a continuous mapping fulfilling for some $p\ge 11/5$, some $0< \underline{\nu}, \overline{\nu}<+\infty$ and for all $\thet\in \R_+$, $\D_1,\D_2 \in \R^{3\times 3}_{\rm sym}$
\begin{subequations}\label{nu}
\begin{align}
(\S^*(\thet, \D_1)-\S^*(\thet, \D_2)): (\D_1-\D_2)&\ge 0,\label{nu1}\\
\S^*(\thet,\D_1):\D_1\ge \underline{\nu}|\D_1|^p - \overline{\nu}, \quad |\S^*(\thet, \D_1)|&\le \overline{\nu}(1+|\D_1|)^{p-1}, \quad \S^*(\thet, 0)=0.\label{nu2}
\end{align}
\end{subequations}
Note that the prototypic relation $\S\sim \nu(\thet) |\D\vv|^{p-2}\D\vv$ falls into the class \eqref{nu}.

The key motivation of the paper is the following stability result:  to show that $(\vv, \thet)\to (0,\thetah)$ in a suitable topology as $t\to \infty$, where $\thetah$ is the unique solution to
\begin{equation}\label{thetahat}
\begin{split}
-\diver(\kappa(\thetah)\nabla\thetah)&=0 \ \mbox{ in $\Omega$},\qquad \thetah=\thet_b \ \mbox{ on $\partial\Omega$}.
\end{split}
\end{equation}
For linear models of the form $\S=\nu(\thet)\D \vv$, such result was already proven in~\cite{DosPruRaj}  provided that the solution $(\vv, \theta)$ is smooth. However, the existence of smooth (or sufficiently regular) solution is not known and the only available results focus only on weak solutions. In three dimensional setting, the existence of a weak solution was firstly proved for $p\ge 11/5$ in \cite{Consiglieri}, see also the related work \cite{FrMaRu10}. Later, for $p\in (9/5,11/5)$ and slightly different boundary conditions, the existence of a weak solution was established in
\cite{BuMaRa09}, see also \cite{MaZa18} for more complicated model, with one proviso. The identity \eqref{i3} was replaced by the inequality, which in terms of the entropy
\begin{equation}\label{entropy2}
\eta = \log\thet
\end{equation}
can be rewritten into the so-called entropy inequality
\begin{equation}
\label{entropy}
\partial_t \eta + \diver(\eta \vv) + \diver\left(\frac{\vc q}{\thet}\right) \geq \frac{1}{\thet} \S:\D\vv + \kappa(\thet)\frac{|\nabla\thet|^2}{\thet^2}.
\end{equation}
Unfortunately, such a result does not allow us to use the procedure developed in \cite{DosPruRaj} and to prove a stability result. Therefore, we need to change the methods and results used in \cite{Consiglieri,BuMaRa09,BulFeiMal} significantly. Thus, our main goal is to prove the existence of a weak solution, which satisfies \eqref{entropy} with the \emph{equality sign} and also to show that the temperature is continuous with respect to time into the topology of $L^1(\Omega)$, which is the natural function space.  Then, inspired by \cite{DosPruRaj}, we know that one can renormalize \eqref{entropy} by a properly chosen set of functions. Indeed, in the standard approach of renormalization, based just on \eqref{i3}, one needs (due to the commutator lemma) that $\thet \in L^{p'}(Q)$. Unfortunately, this is true only for $p>\frac52$. Therefore, we introduce \eqref{entropy} with equality sign and then to renormalize it, one just requires $\eta\in L^{p'}(Q)$, which is true for any $p>1$.  Our result can be understood as a starting point and the corner stone of further rigorous stability results in various models of incompressible heat conducting fluids. Furthermore, it seems that we identify the proper notion of solution for which we are able to show its existence and also, in future, to show the convergence to equilibrium.

%we requireHence
%
%
%
%\begin{equation}\label{eq-t}
%\dert\thet+\diver(\thet\vv) - \diver(k(\thet)\nabla\thet) = \S:\D{\vv}
%\end{equation}
%
%Regarding the internal energy we assume
%\begin{equation}
%e = \cni \thet \mbox{ with }\cni\in (0, \infty)
%\end{equation}
%and as a consequence the entropy satisfies
%
%Finally we assume the heat flux fulfills the Fourier law
%
%
%By virtue of the constitutive relations introduced, equation \eqref{i3} becomes
%
%and we set $\cni=1$.
%In addition to the balance laws \eqref{bl}, one requires that the second law of the thermodynamics is fulfilled:
%
%

\section{Rigorous statement of the main result}

%\subsection{Renormalized entropy solutions}
We use in what follows the standard notation for the Lebesgue, the Sobolev and the Bochner spaces and endow them with standard norms. The symbol $C^{\infty}_0$ is reserved for smooth compactly supported functions and the function spaces related to the incompressible setting are denoted by   $W^{1,p}_{0, \diver}:=\{\vv \in W^{1,p}_{0}(\Omega; \R^3);\; \diver \vv =0\}$ and $L^2_{0,\diver}$ denotes the closure of $W^{1,2}_{0,\diver}$ in $L^2$ topology. Duality pairing between $W^{1,p}_{0, \diver}$ and their duals is denoted $\langle\cdot,\cdot\rangle$.

Next, we prescribe the assumptions on data. Recall, we assume that $\vc{q}$ and $\S^*$ satisfy \eqref{Fourier}--\eqref{nu}. For initial and boundary data we consider
\begin{equation}\label{data}
\vv_0\in L^2_{0,\diver}, \ \thet_0\in L^1(\Omega),  \ \thetah \in W^{1,2}(\Omega)\cap L^{\infty}(\Omega),
\end{equation}
where $\thetah$ is the solution to \eqref{thetahat}. Hence, we transferred all assumptions on the boundary behaviour of $\thet_b$ to the uniquely defined $\thetah$. Finally, we suppose that
\begin{equation}\label{min}
\mu:=\min\left\{\essinf_{x\in \Omega}\thetah(x), \; \essinf_{x\in \Omega}\thet_0(x)\right\}>0.
\end{equation}

The main result of the paper is following.
\begin{theorem}[Existence of a solution fulfilling entropy equality]\label{thm}
Let $\Omega\subset \R^3$ be a bounded domain with Lipschitz boundary. Assume that $\S^*$ and $\kappa$ satisfy \eqref{k}--\eqref{nu} with $p\geq 11/5$. Then for any data $\vv_0, \ \thet_0, \thetah$ fulfilling \eqref{data}--\eqref{min} there exists a quadruplet $(\vv, \S, \thet, \eta)$ fulfilling
\begin{align}
\vv&\in C([0, T]; L^2_{0,\diver})\cap L^p(0, T; W_{0, {\rm div}}^{1, p}),\\
  \partial_t\vv &\in L^{p'}(0, T; (W_{0, {\rm div}}^{1, p})^*), \ \S\in L^{p'}(Q; \R^{3\times 3}),\\
  \thet &\in C([0, T]; L^1(\Omega)), \ (\thet)^\alpha \in  L^2(0, T; W^{1,2}(\Omega)) &&\mbox{ for any } \alpha\in (0, {1}/{2}),\\
 \thet&\in L^r(Q) &&\mbox{ for any } r\in [1, 5/3), \\
\thet - \thetah &\in L^s(0, T; W_0^{1, s}(\Omega))  &&\mbox{ for any } s\in [1, 5/4),\\
\eta &\in L^2(0,T; W^{1,2}(\Omega))\cap L^q(Q) &&\mbox{ for any } q\in [1, +\infty),
\end{align}
and satisfying \eqref{bl} and \eqref{entropy} in the following sense:
\begin{itemize}
\item[] {\bf Momentum equation:} The Cauchy stress is of the form $\S=\S^*(\thet, \D\vv)$ a.e. in $Q$, the initial datum fulfils $\vv(0)=\vv_0$ and for all $\vw\in L^p(0, T; W_{0, {\rm div}}^{1,p})$
\begin{equation}\label{T1}
\begin{aligned}
  \int_0^T{\langle\partial_t\vv, \vw\rangle}\dt + \int_0^T\intO{\S:\D\vw}\dt = \int_0^T\intO{(\vv\otimes\vv): \D\vw}\dt;
  \end{aligned}
\end{equation}

\item[]{\bf Internal energy balance:} Temperature satisfies the minimum principle $\thet\geq \mu \mbox{ a.e. in } Q$, the initial condition fulfils $\thet(0)=\thet_0$ and for all $\varphi\in C^\infty_0((-\infty, T) \times \Omega)$
  \begin{equation}
  \begin{aligned}\label{T2}
-\int_0^T\intO{\thet \partial_t\varphi}\dt - \int_0^T\intO{\thet\vv\cdot\nabla\varphi}\dt + \int_0^T\intO{\kappa(\thet)\nabla\thet\cdot\nabla\varphi}\dt \\
 = \int_0^T\intO{\S:\D\vv\, \varphi}\dt + \intO{\thet_0\varphi(0)}; %\mbox{ for all } \varphi\in C^\infty_0((-\infty, T); C_0^\infty(\Omega))
%\\
%\thet\geq m \mbox{ a.e. in } (0, T) \times\Omega, \ \thet=\thet_b \mbox{ on } (0, T)\times \partial\Omega {\color{blue} FIX how}.
\end{aligned}
\end{equation}

\item[] {\bf Entropy equation:} Entropy is given as $\eta =\ln \thet \mbox{ a.e. in } Q$, $\eta_0:= \ln \thet_0$ and for all $\varphi \in C_0^\infty((-\infty, T)\times\Omega)$
   \begin{equation}\label{entropy-limit}
   \begin{aligned}
-\intTO{\eta\dert\varphi}  - \intTO{\eta\,\vv\cdot \nabla\varphi}  + \intTO{\kappa(\thet)\nabla\eta\cdot\nabla\varphi} \\
= \intTO{\frac{1}{\thet}\,\S:\D{\vv}\,\varphi} + \intTO{\kappa(\thet)\frac{|\nabla\thet|^2}{(\thet)^2}\,\varphi} + \intO{\eta_0\, \varphi(0)}.\end{aligned}
\end{equation}
%\end{itemize}
%\begin{itemize}
%\item[] {\bf Renormalized entropy equation:} {\color{red}TODO $\varphi \in$ initial data} The relative entropy is given as $\etad:= \ln ({\thet}/{\thetah})$ and for any bounded function $\mathcal{T}=\mathcal{T}(s)$ defined for any $s>0$ such that $\mathcal{T}\in C^{0,1}(0,+\infty)$, $\mathcal{T}(1)=0$ and $\T'=0$ on some $(C,+\infty)$ and for $F=F(s)$ such that $F'(s)= \mathcal{T}(s)/s$ for any  $s>0$ it holds
%   \begin{equation}\label{renormalized-entropy}
%\begin{split}
%&- \intTO{\dert\varphi\, \thetah \,F\left(\frac{\thet}{\thetah}\right)} -\intTO{\nabla\eta\cdot\vv \,\varphi \, \thetah\, \mathcal{T}\!\left(\frac{\thet}{\thetah}\right)}\\
%& + \intTO{\kappa(\thet){|\nabla\etad|^2}\varphi\, \thetah\, \left[\mathcal{T}'\left(\frac{\thet}{\thetah}\right)\frac{\thet}{\thetah}-\mathcal{T}\left(\frac{\thet}{\thetah}\right) \right]} \\
%&= \intTO{\frac{1}{\thet}\,\S:\D{\vv}\,\varphi \,\thetah\, \mathcal{T}\left(\frac{\thet}{\thetah}\right)}. %\mbox{ for all } \varphi \in C^\infty_0((0, T)), \ \etad:=\log \frac{\thet}{\thetah}.
%\end{split}
%\end{equation}
\end{itemize}
\end{theorem}

%\subsection{Main theorems}

%{\color{red} Maybe add two selling/explaining sentences}

\section{Proof of Theorem~\ref{thm}}
%\subsection{Existence of approximations}
The existence proof relies on the methods developed in \cite{BuMaRa09} and a large part of the proof is identical. Therefore, we omit unnecessary details and focus mainly on the new aspects of the proof, i.e., on the proof of entropy equality~\eqref{entropy-limit}.

Hence, following \cite{BuMaRa09} (compare also with \cite{Consiglieri}, where a different approach is used), we introduce a basis  $\{ \vw_j\}_{j=1}^\infty$ of the space $W^{3,2}(\Omega; \R^3)\cap  W^{1,p}_{0,\diver}$ that is orthonormal in $L^2_{0,\diver}$, see \cite[Appendix A.4]{MaBook}. Next, for given initial conditions $\vv_0$ and $\thet_0$, we denote $\vv^{n}_0$ the projection of $\vv_0$ onto the subspace $[\vw_1,\dots, \vw_n]$, and $\thet^n_0\in L^2(\Omega)$, fulfilling $\thet^n_0\ge \mu$ a.e. in $\Omega$,  is the regularization of $\thet_0$ such that
\begin{align}
 \vv_0^{n}&\to \vv_0  \ &&\mbox{ strongly in } L^2_{0,\diver} \mbox{ as } n\to +\infty,\\
\thet^n_0 &\to \thet_0\  &&\mbox{ strongly in } L^1(\Omega) \mbox{ as }n\to +\infty.\label{thet-zero}%\\
%&     \thet_0^n(x)\geq m &&\mbox{ for a.a. } x\in \Omega.
\end{align}
Then, for every $n\in \mathbb{N}$, we can find a triple $(\vv^n, \thet^n, \S^n)$, such that
%\begin{theorem}\label{thm:ex-approx}
%For every $n\in\N$ there is
$\vv^n\in W^{1,2}((0,T), W^{3,2}(\Omega;\R^3)\cap W^{1,p}_{0,\diver})$, $\theta^n\in L^\infty((0, T); L^2(\Omega))\cap L^2((0, T); W^{1,2}(\Omega))\cap W^{1,2}((0, T); (W^{1,2}_0(\Omega))^*)$ and $\S^n \in L^{\infty}(Q)$ and such that for a.e. $t\in(0,T)$
\begin{align}
\label{ode11}
&\intO{[\partial_t\vv^{n} \cdot \vw_j  - (\vv^{n} \otimes \vv^{n}): \nabla \vw_j + \S^{n}: \D\vw_j]} = 0 \qquad \mbox{ for all } j=1,\dots,n,
\end{align}
and for all $\psi\in L^2(0,T;W_0^{1,2}(\Omega))$
\begin{align}
\label{ode22}
\begin{aligned}
&\int_0^T\dual{\partial_t\thet^{n}}{\psi}\dt +\intTO{[-\thet^{n} \vv^{n} \cdot \nabla \psi + \kappa(\thet^n)\nabla\thet^{n}\cdot \nabla \psi]} \\
&= \int_0^T\intO{\S^{n}: \D\vv^{n} \psi} \dt.
\end{aligned}
\end{align}
In addition, $\vv^n$ and $\S^n$ are given by
\begin{equation}
  \vv^n= \sum_{i=1}^n c^{n}_i(t)\vw_i(x) \quad \mbox{ and } \quad   \S^n= \S^*(\thet^n, \D \vv^n),
  \end{equation}
the initial conditions $\vv^n(0, \cdot)=\vv^n_0$, $\thet^n(0, \cdot)=\thet_0^n$ are satisfied, $\thet^n$ attains the boundary conditions, i.e., $\thet^n_{|\partial\Omega}=\thetah$ and fulfills the minimum principle, i.e., $\thet^n \ge \mu$ a.e. in $Q$.
%\end{equation}
%\end{theorem}

Then, following e.g. \cite{BuMaRa09} and defining $\thetat^n:=\thet^n - \thetah$, it is rather standard to deduce the following $n$-independent a~priori estimates valid for all $r\in [1,5/3)$, $s\in [1, 5/4)$ and $\alpha \in (0,1/2)$
%
%By standard procedure one can derive following a~priori estimates for approximative solution from Theorem~\ref{thm:ex-approx}.
%??? It is presented in the Appendix.
%??? Are we interested in dependence of constants on $T$? Some os them depend, some not.
%
%\begin{theorem}\label{thm:apr1}
%  Let $\vv^n$ and $\thet^n$ be as in Theorem~\ref{thm:ex-approx}, $\thetat^n=\thet^n-\hat{\thet}$. The next estimates are independent of $n$.
\begin{align}
    \label{apr:vn}
   % \begin{gathered}
      \nm{\vv^{n}}_{L^\infty(0,T; L^2_{0,\diver})}+\nm{\vv^{n}}_{L^p(0,T; W^{1,p}_{0,\diver})}+\nm{\S^n}_{L^{p'}(Q)}+\nm{\vv^{n}}_{L^{\frac{5p}3}(Q)}&\leq C
    %\end{gathered}
    \\
   \label{timen}
 \nm{\thetat^n}_{L^\infty(0,T;L^1(\Omega))}+\nm{\partial_t \thet^n}_{L^1(0, T; (W^{1, 10}_0(\Omega))^*)}+\nm{\partial_t\vv^n}_{L^{p'}(0, T; (W^{1,p}_{0,\diver})^*)} &\leq C,\\
\label{eq:re3}
   \nm{(\theta^n)^\alpha}_{L^2(0,T;W^{1,2}(\Omega))}+\nm{\thetat^n}_{L^r(Q)}+\nm{\thetat^n}_{L^s(0, T; W_0^{1,s}(\Omega))}&\leq C(\alpha, r, s).%\\
%\label{uniform-alpha}
%\nm{(1+|\thetat^n|)^\alpha}_{L^2(0,T;W^{1,2}(\Omega))}&\leq C(\alpha),\\
%\label{uniform-r}
%\nm{\thetat^n}_{L^r((0, T)\times\Omega)}+\nm{\thetat^n}_{L^s(0, T; W_0^{1,s}(\Omega))}&\leq C(r,s)%,\\
%\nm{\partial_t \thet^n}_{L^1((0, T); (W^{10}_0(\Omega))^*)} &\leq C.%\\
%\label{thet-below}
%   \thet^n\geq m \mbox{ a.e. in } (0, T) \times\Omega.
\end{align}
% ??? What is $q$ sufficiently large
%\end{theorem}
%\begin{proof}
%  The estimates are obtained in a standard way. Compare articles \cite{???} for similar estimates in a slightly different setting. The detailed proof is presented in the appendix.
%\end{proof}

This is the starting point of the proof. It is rather sketchy, but it does not contain any essentially new information. However, it is not the case for what follows and therefore we prove it in full details.

\subsection{Limit in momentum and energy equations as \texorpdfstring{$n\to+\infty$}{n}}

By virtue of the established uniform estimates \eqref{apr:vn}--\eqref{eq:re3} and employing the Aubin--Lions compactness lemma, we can extract a subsequence that we do not relabel and we can find $(\vv, \thet, \S)$ such that
\begin{align}
 \vv^n &\rightharpoonup^* \vv && \mbox{ weakly-* in } L^\infty(0, T; L^2_{0,\diver}), \label{Linfty2}\\
 \vv^n &\rightharpoonup \vv &&\mbox{ weakly in } L^p(0, T; W^{1,p}_{0,\diver})\cap W^{1,p'}(0,T; (W^{1,p}_{0,\diver})^*),\label{Lp}\\
 \S^n &\rightharpoonup \S &&\mbox{ weakly in } L^{p'}(Q; \R^3),\label{Lp-S}\\
 %\partial_t\vv^n &\to  \partial_t\vv &&\mbox{ weakly-* in } L^{p'}(0, T; (W_{0, {\rm div}}^{1,p}(\Omega; \R^3))^*),\\
 \vv^n&\to\vv &&\mbox{ strongly in }L^q(Q; \R^3) \mbox{ for any $q\in [1, {5p}/{3})$, and a.e. in $Q$} \label{v-strong}\\
 (\thet^n)^\alpha&\rightharpoonup   (\thet)^\alpha &&\mbox{ weakly in } L^2(0, T; W^{1,2}(\Omega)) \mbox{ for any } \alpha\in (0, {1}/{2}),\label{thet-weak}\\
 \thet^n &\to \thet &&\mbox{ strongly in $L^r(Q)$ for any $r\in [1, 5/3)$, and a.e. in $Q$},\label{thet-strong}\\
  \thet^n -\thetah &\rightharpoonup \thet -\thetah && \mbox{ weakly in } L^s(0, T; W_0^{1, s}(\Omega))  \mbox{ for any } s\in [1, 5/4).\label{grad-thet-weak}
   %\thet^n &\rightharpoonup \thet && \mbox{ weakly in } L^s(0, T; W^{1, s}(\Omega))  \mbox{ for any } s\in [1, 5/4),\label{grad-thet-weak}\\
   %\label{conv-point}
%\vv^n\to \vv &\mbox{ and } \thet^n\to \thet &&\mbox{ a.e. in } (0, T)\times\Omega.
 \end{align}
%??? Do we need all of them
In addition, using \eqref{thet-strong} and the first part of the uniform estimate \eqref{timen}, we have that
\begin{equation}
 \begin{split}\label{theinf1}
  \thet\in L^{\infty}(0,T; L^1(\Omega)).
\end{split}
\end{equation}

Now we are in position to analyze the limit of the formulation \eqref{ode11}-\eqref{ode22}. Indeed,
by virtue of the convergence results \eqref{Lp}-\eqref{v-strong}, we may follow \cite{BleMalRaj} to take the limit in the formulation \eqref{ode11} to deduce that  for all  $\vw\in L^p(0, T; W^{1,p}_{0,\diver})$
\begin{equation}
 \begin{split}\label{AF1}
  \int_0^T{\langle\partial_t\vv, \vw\rangle}\dt + \int_0^T\intO{\S:\D\vw}\dt = \int_0^T\intO{(\vv\otimes\vv): \D\vw}\dt.
\end{split}
\end{equation}
In addition, we have $\vv\in C([0, T]; L^2_{0,\diver})$ and $\vv(0, \cdot)=\vv_0$. To complete the proof of \eqref{T1}, it remains to show
\begin{equation}\label{Minty1}
  \S= \S^*(\thet, \D\vv) \qquad \mbox{ a.e. in } Q.
\end{equation}
Next, thanks to \eqref{v-strong}--\eqref{grad-thet-weak}, one may follow a standard procedure (see e.g. \cite{BulFeiMal}) and let $n\to \infty$ in \eqref{ode22} to show that for all  $\psi\in C^\infty_0((-\infty, T)\times \Omega)$ there holds (compare with \eqref{T2})
\begin{equation}
\begin{split}
-\int_0^T&\intO{\thet \partial_t\psi}\dt - \int_0^T\intO{\thet\vv\cdot\nabla\psi}\dt + \int_0^T\intO{\kappa(\thet)\nabla\thet\cdot\nabla\psi}\dt \\
&= \int_0^T\intO{\S:\D\vv\, \psi}\dt + \intO{\thet_0\psi(0)} \label{3.18}
\end{split}
\end{equation}
provided that we show
\begin{equation}
   \S^n:\D{\vv}^n \rightharpoonup \S:\D{\vv} \mbox{ weakly in } L^1(Q).\label{Sn}
\end{equation}
Obviously, we also have $\thet\geq \mu$ a.e. in $Q$. It remains to show \eqref{Minty1} and \eqref{Sn}. Note that both of them are consequences of the Minty method. Indeed, it follows from \eqref{ode11} that (using the fact that $\diver \vv^n=0$ and \eqref{Lp})
$$
\begin{aligned}
&\limsup_{n\to \infty}\int_0^T \intO{\S^n\cdot \D \vv^n}\dt =\limsup_{n\to \infty} -\int_0^T \langle \partial_t \vv^n, \vv^n \rangle \dt\\
&=\limsup_{n\to \infty}\left( -\int_0^T \langle \partial_t (\vv^n-\vv), \vv^n-\vv \rangle \dt-\int_0^T \langle \partial_t \vv, \vv^n-\vv \rangle \dt-\int_0^T \langle \partial_t \vv^n, \vv \rangle \dt \right)\\
&\le %\limsup_{n\to \infty}\left( \frac12 \|\vv^n(0)-\vv_0\|_2^2 -\int_0^T \langle \partial_t \vv, \vv^n-\vv \rangle \dt-\int_0^T \langle \partial_t \vv^n, \vv \rangle \dt\right) =
\frac12 (\|\vv_0\|_2^2- \|\vv(T)\|_2^2)=\int_0^T \intO{\S : \D \vv} \dt,
\end{aligned}
$$
where the last identity follows from \eqref{AF1} with setting $\vw:=\vv$. Consequently, using this estimate, the monotonicity and the growth assumption \eqref{nu}, the strong convergence result \eqref{thet-strong}, the weak convergence \eqref{Lp} and the Lebesgue dominated convergence theorem, we deduce that for all $\overline{\D}\in L^{p}(Q; \R^{3\times 3})$ there holds
\begin{equation}\label{Mintyc}
\begin{split}
0&\le \limsup_{n\to \infty} \int_0^T \intO{(\S^n-\S^*(\thet^n,\overline{\D})):(\D\vv^n -\overline{\D}) }\dt\\
&\le  \int_0^T \intO{(\S-\S^*(\thet,\overline{\D})):(\D\vv -\overline{\D}) }\dt.
\end{split}
\end{equation}
The classical Minty method then leads to \eqref{Minty1}. Moreover, setting $\overline{\D}:=\D\vv$ in \eqref{Mintyc}, we have
\begin{equation}\label{Mintyd}
\begin{split}
&\limsup_{n\to \infty} \int_0^T \intO{\left|(\S^n-\S^*(\thet^n,\D\vv)):(\D\vv^n -\D\vv)\right| }\dt=0.%\\
%&=\limsup_{n\to \infty} \int_0^T \intO{(\S^n-\S^*(\thet^n,\D\vv)):(\D\vv^n -\D\vv)}\dt=0.
\end{split}
\end{equation}
Hence,
\begin{align}\label{Mintye}
(\S^n-\S^*(\thet^n,\D\vv)):(\D\vv^n -\D\vv) &\to 0 &&\textrm{ strongly in } L^1(Q).
\end{align}
Since
\begin{align*}
\S^*(\thet^n,\D\vv):(\D\vv^n -\D\vv) &\rightharpoonup 0 &&\textrm{ weakly in } L^1(Q),
\end{align*}
which follows from \eqref{Lp}, \eqref{thet-strong} and \eqref{nu2}, we see that \eqref{Mintye} directly implies \eqref{Sn}.

\subsection{Limit in entropy equation as \texorpdfstring{$n\to+\infty$}{n}}

In this section we show the validity of \eqref{entropy-limit}. To see this, we first set $\psi:=\varphi/\thet^n$ in  \eqref{ode22} with arbitrary  $\varphi \in C_0^\infty((-\infty, T)\times\Omega)$  to derive the following identity for approximated entropy $\eta^n:=\ln \thet^n$
\begin{equation}
\begin{split}\label{entropy-weak}
&-\intTO{\eta^n\dert\varphi+\eta^n\,\vv^n\cdot \nabla\varphi-\kappa(\thet^n)\nabla\eta^n \cdot\nabla\varphi}\\
&= \intTO{\frac{\S^n:\D{\vv}^n}{\thet^n}\,\varphi+\kappa(\thet^n)\frac{|\nabla\thet^n|^2}{(\thet^n)^2}\,\varphi} +\intO{\eta^n_0\, \varphi(0)},
\end{split}
\end{equation}
where we set $\eta^n_0:=\ln \thet^n_0$. Next, we want to let $n\to \infty$. The identification of the limit in the terms on the left hand side is rather standard and follows from the convergence results \eqref{Lp}, \eqref{v-strong}, \eqref{thet-strong} and \eqref{grad-thet-weak}. Similarly, to pass to the limit in the first term on the right hand side of \eqref{entropy-weak} is straightforward thanks to \eqref{thet-strong}, \eqref{Sn}, the Egorov and Dunford-Pettis theorems. Also the limit passage in the last term is obvious.  The most problematic term is however the second term on the right hand side since $\kappa(\thet^n)|\nabla\thet^n|^2/(\thet^n)^2$ is uniformly bounded only in $L^1((0, T)\times\Omega)$ and so we cannot even a~priori extract an $L^1$ weakly convergent subsequence. We overcome this in two steps. First, the point-wise convergence of $\nabla\thet^n$ is shown and then the strong convergence of $\kappa(\thet^n)|\nabla\eta^n|^2$ in $L^1((0, T)\times\Omega)$ is deduced.

\subsubsection{Almost everywhere convergence of \texorpdfstring{$\nabla\thet^n$}{n}}
\newcommand{\tnk}{\thet^n_K}
\newcommand{\tmk}{\thet^m_K}
\newcommand{\wmn}{w^{m,n}}
\newcommand{\F}{\mathcal F}
\newcommand{\G}{\mathcal G}
\newcommand{\wmnd}{w^{m,n}_{\delta}}
We start this part with definition of auxiliary cut-off functions. For arbitrary  $k >0$, we define
\begin{equation}\label{def:tk}
\mathcal{T}_k(z) := {\rm sign}(z) \min\{|z|,k\}.
\end{equation}
Its primitive function attaining zero at zero is denoted $\mathcal G_k$, i.e. $\mathcal G_k'=\mathcal{T}_k$, $\mathcal G_k(0)=0$. Note that $|\G_k(s)|\leq k|s|$ for all $s\in\R$. Next, we also introduce a mollification of $\mathcal{T}_k$. For arbitrary  $\delta\in (0, k)$ (typically $\delta \ll 1$), we denote by $\mathcal{T}_{k, \delta}\in C^2(\mathbb{R})$ a mollification of $\mathcal T_k$, which is given by a convolution with a symmetric, positive kernel of radius $\delta$. Such a mollification then has the following properties
\begin{align*}
&\mathcal{T}_{k, \delta}(z)= \mathcal{T}_k(z) \qquad \mbox{ if }   |z|\leq k-\delta \mbox{ or } |z|\geq k+\delta, \qquad | \mathcal{T}^{\,''}_{k, \delta}|\leq C\delta^{-1}\\
  &0\leq \mathcal{T}^{\,'}_{k, \delta}\leq 1, \quad \mathcal{T}^{\,''}_{k, \delta}\leq 0, \quad \mathcal{T}_{k, \delta}\leq \mathcal{T}_k\quad \mbox{on $(0,+\infty)$.}
%\  \mathcal{T}^{\,'}_{k, \delta}(z)=0 \mbox{ for } 0\leq z\leq k-\delta \mbox{ or } |z|\geq k+\delta.
\end{align*}

%Let $K>0$, $\alpha\in(41/66,1)$ (the lower bound comes from considerations above \eqref{lim}) and $G:(0,+\infty)\to(0,+\infty)$ be a function with following properties: $0\leq G(s)\leq C s^{1-\alpha}$, $G$ smooth, increasing, concave, $G'(s)\leq Cs^{-\alpha}$ $(m,+\infty)$ and $|G''(s)|\leq Cs^{-1-\alpha}$ for $s\in(m,+\infty)$.

We fix $m,n,k\in\N$, $k>2\nm{\thetah}_{L^\infty(\partial\Omega)}$ and $\varepsilon, \delta>0$, $\varepsilon<k$ and define $ w^{m,n}_{\delta}=\mathcal{T}_{k+\varepsilon, \delta}(\thet^n)-\mathcal{T}_{k, \delta}(\thet^m)$. We set  $\psi:=\mathcal{T}^{\,'}_{k+\varepsilon, \delta}(\thet^n)\mathcal{T}_\ep(w^{m, n}_{\delta})$ in \eqref{ode22} for $\thet^n$ and $\psi:=\mathcal{T}^{\,'}_{k, \delta}(\thet^m)\mathcal{T}_\ep(w^{m, n}_{\delta})$  in \eqref{ode22} for $\thet^m$. Note that it is allowed since $\mathcal{T}^{\,'}_{k+\varepsilon, \delta}(\thet^n)\mathcal{T}_\ep(w^{m, n}_{\delta})$ and also $\mathcal{T}^{\,'}_{k, \delta}(\thet^m)\mathcal{T}_\ep(w^{m, n}_{\delta})$ belong to $L^2((0, T); W_0^{1,2}(\Omega))$. Then, we subtract the so obtained equations to get
\begin{equation}\label{Tep}
\begin{split}
&\intTO{\kappa(\thet^n)\nabla(\wmnd)\cdot \nabla \mathcal{T}_\varepsilon(\wmnd)}\\
=&\intTO{\bigl(\kappa(\thet^n)-\kappa(\thet^m)\bigr)\nabla \mathcal{T}_{k, \delta}(\thet^m)\cdot\nabla\mathcal{T}_\varepsilon(\wmnd)} -\int_0^T  \langle \dert w^{m,n}_{\delta},  \mathcal{T}_\varepsilon(w^{m,n}_{\delta})\rangle \dt\\
&\phantom{=} + \intTO{ G^{m,n} \mathcal{T}_\varepsilon(\wmnd)+%}  \intTO{
\left[\mathcal{T}_{k+\varepsilon, \delta}(\thet^n) \vv^n - \mathcal{T}_{k, \delta}(\thet^m) \vv^m\right] \cdot \nabla \mathcal{T}_\varepsilon(w^{m,n}_\delta)},
\end{split}
\end{equation}
where we denoted
\begin{equation*}
\begin{split} G^{m,n}:&=\left[\mathcal{T}^{\,'}_{k+\varepsilon, \delta}(\thet^n) \,\S^n:\D\vv^n - \mathcal{T}^{\,'}_{k, \delta}(\thet^m) \,\S^m:\D\vv^m\right]\\
  &-\left[ \kappa(\thet^n) |\nabla \thet^n|^2 \mathcal{T}^{\,''}_{k+\varepsilon, \delta}(\thet^n)- \kappa(\thet^m)|\nabla \thet^m|^2 \mathcal{T}^{\,''}_{k, \delta}(\thet^m)\right].
\end{split}
\end{equation*}
Our first goal is to let $\delta \to 0_+$ in \eqref{Tep}. Note that such convergence procedure is very standard in all terms except the term $G^{m,n}$ involving the second derivative of  $\mathcal{T}_{\cdot,\delta}$. To get a proper $\delta$-independent bound, we consider $M\ge \|\thetah\|_{L^{\infty}(\partial \Omega)}$ and $\delta \in (0,M/2)$. Since $\thet^n=\thetah$ on $\partial \Omega$, we can deduce that $\psi:=1-\mathcal{T}^{\,'}_{M, \delta}(\thet^n)\in L^2(0,T;W^{1,2}_0(\Omega))$ and therefore it can be used in \eqref{ode22}. Such choice then leads to
\begin{equation}
\begin{split}\label{T22}
&\intTO{ \kappa(\thet^n)|\mathcal{T}^{\,''}_{M, \delta}(\thet^n)| |\nabla\thet^n|^2} =\intO{(\thet^n_0-\mathcal{T}_{M, \delta}(\thet^n_0))} \\
&-\intO{(\thet^n-\mathcal{T}_{M, \delta}(\thet^n))(T)} + \intTO{(1-\mathcal{T}^{\,'}_{M, \delta}(\thet^n)) \S^n:\D{\vv}^n }\\
&\leq \int_{\{\thet^n_0>M/2\}}\!\!\!\!\! \thet^n_0\dx + \int_{\{|\thet^n|>M/2\}}\!\!\!\!\!\S^n:\D{\vv}^n \dx\dt\le C,
\end{split}
\end{equation}
where we exploited that $\diver \vv^n=0$, used  the properties of $\mathcal{T}_{M, \delta}$ (in particular the concavity) and \eqref{apr:vn} and \eqref{thet-zero}. In a very similar manner we can estimate the term with time derivative in \eqref{Tep}. Since $\thet^n$ and $\thet^m$ belong for fix $n,m$ to $C([0,T]; L^2(\Omega))$, we have (using the nonnegativity of $\G_{\varepsilon}$ as well as the estimate $\G_{\varepsilon}(s)\le s\varepsilon$)
\begin{equation}\label{dert}
\begin{split}
-\int_0^T \langle \dert w^{m,n}_{\delta},  \mathcal{T}_\varepsilon(w^{m,n}_{\delta})\rangle\dt  =
   - \intTO{\dert \G_\varepsilon(w^{m,n}_{\delta})}\le \intO{\G_\varepsilon(w^{m,n}_{\delta}(0))}\leq C\varepsilon.
\end{split}
\end{equation}
%\begin{equation}\label{T22}
%\intTO{k(\thet^n) |\mathcal{T}^{\,''}_{k, \delta}(\thet^n)| |\nabla\thet^n|^2} \leq \int_{\{\thet^n_0>k/2\}}\!\!\!\!\! \thet^n_0\dx + %\int_{\{|\thet^n|>k/2\}}\!\!\!\!\!\S^n:\D{\vv}^n \dx\dt\le C,
%\end{equation}
%where the last inequality is the consequence of \eqref{apr:vn} and \eqref{thet-zero}.
Hence, we can apply estimates \eqref{T22} and \eqref{dert} and then let $\delta \to 0_+$ in the remaining terms of \eqref{Tep} to deduce (recall $k\ge 2\|\thetah\|_{L^{\infty}(\partial \Omega)}$)
\begin{equation}\label{Tep2}
\begin{split}
&\underline{\kappa}\intTO{|\nabla\mathcal{T}_\varepsilon(\wmn)|^2}\le \intTO{\bigl(\kappa(\thet^n)-\kappa(\thet^m)\bigr)\nabla \mathcal{T}_{k}(\thet^m)\cdot \nabla\mathcal{T}_\varepsilon(\wmn)}\\
%&+C\varepsilon\\% + \intTO{%}  \intTO{
%\left[\mathcal{T}_{k, \delta}(\thet^n) \vv^n - \mathcal{T}_{k, \delta}(\thet^m) \vv^m\right] \cdot \nabla \mathcal{T}_\varepsilon(w^{m,n})}\\
&\phantom{=} + \intTO{%}  \intTO{
\left[\mathcal{T}_{k+\varepsilon}(\thet^n) \vv^n - \mathcal{T}_{k}(\thet^m) \vv^m\right] \cdot \nabla \mathcal{T}_\varepsilon(w^{m,n})}+C\varepsilon,
\end{split}
\end{equation}
where $\wmn:=\mathcal{T}_{k+\varepsilon}(\thet^n) - \mathcal{T}_k(\thet^m)$. Next goal is to let $n,m \to \infty$ and finally $\varepsilon \to 0_+$. We start with the second term on the right hand side. Defining $w_{\varepsilon}:=\mathcal{T}_{k+\varepsilon}(\thet) - \mathcal{T}_k(\thet)$ and using \eqref{v-strong}--\eqref{grad-thet-weak}, we deduce
\begin{equation}\label{lim}
\begin{aligned}
\lim_{\varepsilon \to 0_+} &\limsup_{n\to \infty} \limsup_{m\to \infty}  \intTO{\left[\mathcal{T}_{k+\varepsilon}(\thet^n) \vv^n - \mathcal{T}_{k}(\thet^m) \vv^m\right] \cdot \nabla \mathcal{T}_\varepsilon(w^{m,n})}\\
&=\lim_{\varepsilon \to 0_+} \intTO{\left[\mathcal{T}_{k+\varepsilon}(\thet) \vv - \mathcal{T}_{k}(\thet) \vv\right] \cdot \nabla \mathcal{T}_\varepsilon(w_{\varepsilon})}=0.
\end{aligned}
\end{equation}
The remaining term in \eqref{Tep2} is estimated with the help of \eqref{eq:re3} and the H\"{o}lder inequality as follows (note that $\{|w^{m,n}|<\epsilon\}\cap\{\thet^m\leq k\}=\{|\thet^n-\thet^m|<\epsilon\}\cap\{\thet^m\leq k\}$)
\begin{equation}
\begin{split}
&\intTO{\bigl(\kappa(\thet^n)-\kappa(\thet^m)\bigr)\nabla \mathcal{T}_{k}(\thet^m)\cdot \nabla\mathcal{T}_\varepsilon(\wmn)}\\
&\leq 2\intTO{|\kappa(\thet^n)-\kappa(\thet^m)|(|\nabla \thet^n|^2 + |\nabla \thet^m|^2)\chi_{\{|\wmn|<\ep\}\cap \{\thet^m\le k\}}}\\
&\le C(k)\||\kappa(\thet^n)-\kappa(\thet^m)|\chi_{\{|\wmn|<\ep\}\cap \{\thet^m\le k\}}\|_{L^{\infty}(Q)}\\
&\le C(k)\||\kappa(\thet^n)-\kappa(\thet^m)|\chi_{\{|\thet^n-\thet^m|<\ep\}\cap \{\thet^m\le k\}}\|_{L^{\infty}(Q)}\\
&\le C(k)\sup_{l,s\in [\mu,2k]; \, |s-l|\le \varepsilon}|\kappa(l)-\kappa(s)| \to 0
\end{split}\label{small}
\end{equation}
as $\varepsilon \to 0_+$ thanks to the uniform continuity of $\kappa$ on $[\mu,2k]$.

Hence, using \eqref{lim}--\eqref{small} in \eqref{Tep2}, have
\begin{equation*}%\label{Tep3}
\begin{split}
\lim_{\varepsilon \to 0_+} \limsup_{n\to \infty} \limsup_{m\to \infty} \intTO{|\nabla\mathcal{T}_\varepsilon(\wmn)|^2}=0,
\end{split}
\end{equation*}
which, thanks to the weak lower semicontinuity of the norm in $L^2(Q)$, gives
\begin{equation}\label{Tep3}
\begin{split}
\lim_{\varepsilon \to 0_+} \limsup_{n\to \infty}  \intTO{|\nabla\mathcal{T}_\varepsilon(\mathcal{T}_{k+\varepsilon}(\thet^n)-\mathcal{T}_k(\thet))|^2}=0.
\end{split}
\end{equation}
Finally, we use the H\"{o}lder inequality, the uniform bound \eqref{eq:re3} and the Tschebyschev inequality to get
$$
\begin{aligned}
&\intTO{|\nabla \thet^n -\nabla \thet|}\le \intTO{|\nabla \mathcal{T}_{\varepsilon}(\thet^n -\thet)|}+ \intTO{|\nabla \thet^n -\nabla \thet|\chi_{\{|\thet^n-\thet|>\varepsilon\}}}\\
&\le \intTO{|\nabla \mathcal{T}_{\varepsilon}(\mathcal{T}_{k+\varepsilon}(\thet^n)-\mathcal{T}_k(\thet))|}+\intTO{|\nabla \thet^n -\nabla \thet|\chi_{\{|\thet^n-\thet|>\varepsilon\}\cup \{|\thet^n|+|\thet|>k/2\}}}\\
&\le C\left(\intTO{|\nabla \mathcal{T}_{\varepsilon}(\mathcal{T}_{k+\varepsilon}(\thet^n)-\mathcal{T}_k(\thet))|^2}\right)^{\frac12}\\
&\quad +C\|\nabla \thet^n -\nabla \thet\|_{L^{\frac98}(Q)}\left(|\{|\thet^n-\thet|>\varepsilon\}|^{\frac19}+|\{|\thet^n|+|\thet|>k/2\}|^{\frac19}\right)\\
&\le C\left(\intTO{|\nabla \mathcal{T}_{\varepsilon}(\mathcal{T}_{k+\varepsilon}(\thet^n)-\mathcal{T}_k(\thet))|^2}\right)^{\frac12}+\frac{C\|\thet^n-\thet\|^{\frac19}_{L^1(Q)}}{\varepsilon^{\frac{1}{9}}}+\frac{C}{k^{\frac{1}{9}}}.
\end{aligned}
$$
Therefore, the convergence result \eqref{thet-strong} leads to
$$
\begin{aligned}
&\lim_{n\to \infty}\intTO{|\nabla \thet^n -\nabla \thet|}\le  C\left(\limsup_{n\to \infty}\intTO{|\nabla \mathcal{T}_{\varepsilon}(\mathcal{T}_{k+\varepsilon}(\thet^n)-\mathcal{T}_k(\thet))|^2}\right)^{\frac12}+\frac{C}{k^{\frac{1}{9}}}.
\end{aligned}
$$
As the left hand side is independent of $\varepsilon$ and $k$, we may let first $\varepsilon \to 0_+$ and use \eqref{Tep3} to eliminate the first term on the right hand side and then let $k\to\infty$ to handle the second term on the right hand side and thus to observe that
\begin{equation}
\nabla \thet^n \to \nabla \thet \qquad \textrm{strongly in } L^1(Q)\label{NS1}
\end{equation}
and consequently (for a subsequence)
%as consequence
%\begin{equation*}
%\nabla (\thet^n)^{1-\alpha} \to \nabla (\thet)^{1-\alpha} \mbox{ a.e. in } (0, T)\times\Omega,
%\end{equation*}
%which together with the second part of \eqref{conv-point} yields
\begin{equation}\label{ae-nabla}
\nabla \thet^n \to \nabla \thet \qquad \mbox{ a.e. in } Q.
\end{equation}
\let\tktn\relax

\subsubsection{Strong convergence of \texorpdfstring{${\kappa(\thet^n)|\nabla \thet^n|^2}/{(\thet^n)^2}$}{n} in \texorpdfstring{$L^1$}{L1}-norm}

We start the proof of the claim by showing a strong convergence of $\nabla \mathcal{T}_k(\thet^n-\thetah)$ in $L^2(Q)$ for arbitrary $k$. To prove such a result, we want to set $\psi:=\mathcal{T}_k(\thet^n-\thetah)$ in \eqref{ode22}, let $n \to \infty$ and compare the limit with \eqref{T2} tested by $\varphi:= \mathcal{T}_k(\thet-\thetah)$. However, such test functions are not allowed in general and therefore we must proceed more carefully. We fix an arbitrary $T^*\in (0,T)$, which is the Lebesgue point of $\thet(t)$ as a function in $L^1(0,T;L^1(\Omega))$.
%
%We recall that the function $\mathcal{G}_k$ denotes the primitive function to $\mathcal{T}_k$, see below \eqref{def:tk}.
%First we realize that at almost every time level
Using the fact that $\diver \vv^n=0$, we see that
\begin{equation}\label{thetab-n}
\begin{split}
    &\intO{\nabla \thet^n \cdot \vv^n \, \mathcal{T}_k(\thet^n -\thetah)}% =
%\intO{\nabla (\thet^n-\thetah) \cdot \vv^n \, \mathcal{T}_k(\thet^n -\thetah)} +
%\intO{\nabla \thetah \cdot \vv^n \, \mathcal{T}_k(\thet^n -\thetah)}\\
%&=\intO{\nabla \mathcal{G}_k(\thet^n-\thetah) \cdot \vv^n} +
%\intO{\nabla \thetah \cdot \vv^n \, \mathcal{T}_k(\thet^n -\thetah)}
=\intO{\nabla \thetah \cdot \vv^n \, \mathcal{T}_k(\thet^n -\thetah)}.
\end{split}
\end{equation}
Then, we set $\psi:=\mathcal{T}_k(\thet^n-\thetah) \chi_{[0,\tau]}$ with $\tau\in (T^*,T)$ arbitrary in \eqref{ode22}. We deduce (since $\thet^n \in C([0,T]; L^2(\Omega))$ and $\thetah$ is independent of time) by using \eqref{thetab-n} that
\begin{equation}\label{MBa}
\begin{split}
&\int_0^{\tau }\intO{\kappa(\thet^n)\nabla (\thet^n-\thetah) \cdot \nabla [\mathcal{T}_k(\thet^n -\thetah)]}\dt= - \intO{\mathcal{G}_k(\thet^n(\tau) -\thetah)-\mathcal{G}_k(\thet^n_0 -\thetah)}\\
& + \int_0^{\tau}\intO{-\nabla \thetah \cdot \vv^n \, \mathcal{T}_k(\thet^n -\thetah)-\kappa(\thet^n)\nabla \thetah \cdot \nabla [\mathcal{T}_k(\thet^n -\thetah)]+\mathcal{T}_k(\thet^n -\thetah)\, \S^n:\D{\vv^n} } \dt.
\end{split}
\end{equation}
Then, it follows from \eqref{v-strong}, \eqref{thet-strong}, \eqref{thet-zero}, \eqref{Sn}  and \eqref{MBa} that for $\delta\in(0,T-T^*)$
\begin{equation}\label{MB}
\begin{split}
&\limsup_{n\to \infty}\int_0^{T^*}\intO{\kappa(\thet^n)\nabla \thet^n \cdot \nabla [\mathcal{T}_k(\thet^n -\thetah)]}\dt\\
&=\limsup_{n\to \infty}\int_0^{T^*}\intO{\kappa(\thet^n)\nabla [\thet^n-\thetah] \cdot \nabla [\mathcal{T}_k(\thet^n -\thetah)]+\kappa(\thet^n)\nabla \thetah \cdot \nabla [\mathcal{T}_k(\thet^n -\thetah)]}\dt\\
&\le \limsup_{n\to \infty} \fint_{T^*}^{T^*+\delta}\int_0^{\tau}\intO{\kappa(\thet^n)\nabla [\thet^n-\thetah] \cdot \nabla [\mathcal{T}_k(\thet^n -\thetah)]}\dt\dtau\\
&\qquad\qquad  +\int_0^{T^*}\intO{\kappa(\thet)\nabla \thetah \cdot \nabla [\mathcal{T}_k(\thet -\thetah)]}\dt\\
%&= - \delta^{-1}\int_{T^*}^{T^*+\delta}\intO{\mathcal{G}_k(\thet(\tau) -\thetah)-\mathcal{G}_k(\thet_0 -\thetah)}\dtau\\
%& + \delta^{-1}\int_0^{\delta}\int_0^{T^*+\tau}\intO{\nabla \thetah \cdot \vv \, \mathcal{T}_k(\thet -\thetah)+k(\thet)\nabla \thet \cdot \nabla [\mathcal{T}_k(\thet -\thetah)]+\mathcal{T}_k(\thet -\thetah)\, \S:\D{\vv} } \dt\dtau\\
&\overset{\eqref{MBa}}{\to} -\intO{\mathcal{G}_k(\thet(T^*) -\thetah)-\mathcal{G}_k(\thet_0 -\thetah)}
 %\\&
+ \int_0^{T^*}\intO{(\S:\D{\vv}-\nabla \thetah \cdot \vv)\mathcal{T}_k(\thet -\thetah) } \dt
  % +\kappa(\thet)\nabla \thet \cdot \nabla [\mathcal{T}_k(\thet -\thetah)]
\end{split}
\end{equation}
as $\delta \to 0_+$.
% that We note that we can extract a subsequence (denoted again by the same symbol) of $\{\thet^n\}_n$ in such a way that $\thet^n(T^*) \to \thet(T^*)$ in $L^1(\Omega)$ for a.e. $T^*\in(0,T)$, see \eqref{thet-strong}. We fix such $T^*$, test \eqref{ode22} by $\mathcal{T}_k(\thet^n -\thetah)\chi_{(0,T^*)}$ and finally take the limit as $n\to+\infty$. We get, using \eqref{v-strong}, \eqref{conv-point}, \eqref{Sn},
% \begin{equation}\label{Gk}\begin{split}
% \Dt\intO{\mathcal{G}_k(\thet^n -\thetah)} + \intO{\nabla \thet^n \cdot \vv^n \, \mathcal{T}_k(\thet^n -\thetah)} &\\
% +\intO{k(\thet^n)\nabla \thet^n \cdot \nabla [\mathcal{T}_k(\thet^n -\thetah)]} = \intO{\mathcal{T}_k(\thet^n -\thetah)\, \S^n:\D{\vv}^n },&
% \end{split}\end{equation}
% but observe that
% Let $T^*$ be such that $\thet^n(T^*) \to \thet(T^*)$ in $L^1(\Omega)$, and let us
% integrate \eqref{Gk} over $(0, T^*)$ and take the limit as $n\to+\infty$, it follows
%\begin{equation}\label{MB}\begin{split}
%&\intO{\mathcal{G}_k(\thet(T^*) -\thetah)} -\intO{\mathcal{G}_k(\thet_0 -\thetah)} + \int_0^{T^*}\intO{\nabla \thetah \cdot \vv \, %\mathcal{T}_k(\thet -\thetah)}\dt\\
%& -\int_0^{T^*} \intO{\mathcal{T}_k(\thet -\thetah)\, \S:\D{\vv} } \dt= - \lim_{n\to +\infty} \int_0^{T^*}\intO{k(\thet^n)\nabla \thet^n \cdot %\nabla [\mathcal{T}_k(\thet^n -\thetah)]}\dt
%\end{split}\end{equation}
%??? Shall we discuss more convergence $\intO{\mathcal{G}_k(\thet^n(T^*) -\thetah)}$ to $\intO{\mathcal{G}_k(\thet(T^*) -\thetah)}$ and the term %with time derivative generally
 Now, the aim is to identify the final expression in \eqref{MB}. We want to set $\psi:=\mathcal{T}_k(\thet-\thetah)$ in \eqref{3.18}. Unfortunately, it is not possible due to low regularity of $\thet$. Therefore, we must proceed differently. In what follows, we always consider $M,k,\ep,\delta>0$ such that $M>k+1+\|\thetah\|_{\infty}$ and  $\delta\in(0,1)$. Due to the smoothness of $\thet^n$, we can set $\psi(\tau,x):=\mathcal{T}'_{M, \delta}(\thet^n(\tau,x)) \chi_{(t,t+\varepsilon)}(\tau)\varphi(x)$, $\tau\in(0,T)$, $x\in\Omega$ in \eqref{ode22} to deduce that for all $\varphi\in W^{1,2}_0(\Omega) \cap L^{\infty}(\Omega)$ and all $t\in (0,T^*)$, $\varepsilon\in(0,T-T^*)$ we have
\begin{equation*}
\begin{split}
    &\intO{\partial_t \fint_t^{t +\ep} \mathcal{T}_{M, \delta}(\thet^n) \dtau \varphi}
    + \intO{\fint_t^{t +\ep} \nabla [\mathcal{T}_{M, \delta}(\thet^n)]\cdot \vv^n \dtau \varphi}\\
    &+ \intO{  \fint_t^{t +\ep} \kappa(\thet^n)\nabla \mathcal{T}_{M, \delta}(\thet^n)\dtau \cdot \nabla \varphi}
     = \intO{\fint_t^{t +\ep}  \mathcal{T}^{\,'}_{M, \delta}(\thet^n) \, \S^n:\D{\vv}^n\dtau \varphi}\\
    &+  \intO{\fint_t^{t +\ep} \kappa(\thet^n) \mathcal{T}^{\,''}_{M, \delta}(\thet^n) |\nabla\thet^n|^2\dtau\varphi}.
    % \ \mathcal{T}_k \left( \fint_t^{t +\ep} \mathcal{T}_{M, \delta}(\thet^n) \dtau - \thetah\right)}\dt
\end{split}
\end{equation*}
Then, we fix a measurable set $J\subset(0,T^*)$ and integrate over $t\in J$. Using the convergence results \eqref{v-strong}--\eqref{thet-strong} and \eqref{Sn} next the density of step functions in $L^2(0,T^*;W^{1,2}_0(\Omega))\cap L^\infty(0,T;L^\infty(\Omega))$, we obtain that for all $\varphi\in L^{2}(0,T; W^{1,2}_0(\Omega))\cap L^{\infty}(Q)$,
\begin{equation}\label{eq:dvakrize}
\begin{split}
    &\int_0^{T^*}\intO{\partial_t \thet^{M,\delta}_{\ep} \varphi}\dt
    + \int_0^{T^*}\intO{\fint_t^{t +\ep} \nabla [\mathcal{T}_{M, \delta}(\thet)]\cdot \vv \dtau \varphi}\dt\\
    &+ \int_0^{T^*}\intO{  \fint_t^{t +\ep} \kappa(\thet)\nabla \mathcal{T}_{M, \delta}(\thet)\dtau \cdot \nabla \varphi}\dt
     = \int_0^{T^*}\intO{\fint_t^{t +\ep}  \mathcal{T}^{\,'}_{M, \delta}(\thet) \, \S:\D{\vv}\dtau \varphi}\dt\\
    &+ \limsup_{n\to \infty} \int_0^{T^*}\intO{\fint_t^{t +\ep} \kappa(\thet^n) \mathcal{T}^{\,''}_{M, \delta}(\thet^n) |\nabla\thet^n|^2\dtau\varphi}\dt,
    % \ \mathcal{T}_k \left( \fint_t^{t +\ep} \mathcal{T}_{M, \delta}(\thet^n) \dtau - \thetah\right)}\dt
\end{split}
\end{equation}
where we denoted $\thet^{M,\delta}_{\ep}(t,x):=\fint_t^{t +\ep} \mathcal{T}_{M, \delta}(\thet(\tau, x)) \dtau$. Note that for every $\ep>0$ the function $\thet^{M,\delta}_{\ep}$ belongs to $W^{1,\infty}(0,T; L^{\infty}(\Omega))$ and so the integral with time derivative is well defined.
From \eqref{T22}, \eqref{thet-zero}, \eqref{thet-strong} and \eqref{Sn} we get
\begin{equation}\label{eq:nextest}
\begin{split}
&\limsup_{n\to\infty}\intTO{ \kappa(\thet^n)|\mathcal{T}^{\,''}_{M, \delta}(\thet^n)| |\nabla\thet^n|^2} \leq \intO{(\thet_0-\mathcal{T}_{M, \delta}(\thet_0))} \\
&\quad + \intTO{(1-\mathcal{T}^{\,'}_{M, \delta}(\thet)) \S:\D{\vv} }
\leq \int_{\{\thet_0>M/2\}}\!\!\!\!\! \thet_0\dx + \int_{\{|\thet|>M/2\}}\!\!\!\!\!\S:\D{\vv} \dx\dt.
\end{split}
\end{equation}

Setting  $\varphi=\mathcal{T}_k( \thet^{M,\delta}_{\ep} - \thetah )$ in \eqref{eq:dvakrize} and using \eqref{eq:nextest} we obtain
\begin{equation*}
\begin{split}
&\int_0^{T^*}\intO{  \fint_t^{t +\ep} \kappa(\thet)\nabla \mathcal{T}_{M, \delta}(\thet)\dtau \cdot \nabla \mathcal{T}_k \left(\thet^{M,\delta}_{\ep} - \thetah\right)}\dt\\
&\ge-\intO{\mathcal{G}_k\left( \thet^{M,\delta}_{\ep}(T^*)- \thetah\right)-\mathcal{G}_k\left( \thet^{M,\delta}_{\ep}(0)- \thetah\right)}\\
& + \int_0^{T^*}\intO{\left(\fint_t^{t +\ep}  \mathcal{T}^{\,'}_{M, \delta}(\thet) \, \S:\D{\vv}-\nabla [\mathcal{T}_{M, \delta}(\thet)]\cdot \vv \dtau  \right) \mathcal{T}_k \left( \thet^{M,\delta}_{\ep} - \thetah\right)}\dt\\
&-k\left(\int_{\{\thet_0>M/2\}}\!\!\!\!\! \thet_0\dx + \int_{\{\thet>M/2\}}\!\!\!\!\!\S:\D{\vv} \dx\dt\right).
\end{split}
\end{equation*}
%and note that the first integral at the left-hand side can be written as
%$$\intO{\mathcal{G}_k\left[\fint_{T^*}^{T^* + \ep} \mathcal{T}_{M, \delta}(\thet^n) \dtau - \thetah \right]} -\intO{\mathcal{G}_k\left[ \fint_0^\ep \mathcal{T}_{M, \delta}(\thet^n) \dtau - \thetah\right]}. $$
%Thanks to the properties of $ \mathcal{T}_{M, \delta}$ and $\mathcal{T}_k$  the limit as $n\to +\infty$ can be easily identified, thus it holds
%\begin{equation}\begin{split}
%&\intO{\mathcal{G}_k\left[\fint_{T^*}^{T^* + \ep} \mathcal{T}_{M, \delta}(\thet) \dtau - \thetah \right]} -\intO{\mathcal{G}_k\left[ \fint_0^\ep \mathcal{T}_{M, \delta}(\thet) \dtau - \thetah\right]}
% \\&
%+ \int_0^{T^*}\intO{\fint_t^{t +\ep} \nabla [\mathcal{T}_{M, \delta}(\thet)]\cdot \vv \dtau\ \mathcal{T}_k \left( \fint_t^{t +\ep} \mathcal{T}_{M, \delta}(\thet) \dtau - \thetah\right)}\dt\\
%& + \int_0^{T^*}\intO{  \fint_t^{t +\ep} k(\thet)\nabla\mathcal{T}_{M, \delta}(\thet)\dtau \cdot \nabla \mathcal{T}_k \left( \fint_t^{t +\ep} \mathcal{T}_{M, \delta}(\thet) \dtau - \thetah\right)}\dt\\
%& = \int_0^{T^*}\intO{\fint_t^{t +\ep}  \mathcal{T}^{\,'}_{M, \delta}(\thet) \, \S:\D{\vv}\dtau \ \mathcal{T}_k \left( \fint_t^{t +\ep} \mathcal{T}_{M, \delta}(\thet) \dtau - \thetah\right)}\dt\\
%&+\lim_{n\to+\infty}  \int_0^{T^*}\intO{\fint_t^{t +\ep} k(\thet^n) \mathcal{T}^{\,''}_{M, \delta}(\thet^n) |\nabla\thet^n|^2\dtau\ \mathcal{T}_k \left( \fint_t^{t +\ep} \mathcal{T}_{M, \delta}(\thet^n) \dtau - \thetah\right)}\dt.
%\end{split}\end{equation}
Since $T^*$ is the Lebesgue point of $\thet$ and since the initial condition is attained (see e.g. \cite{BuMaRa09, Consiglieri}) we can let $\ep\to 0_+$ and obtain
\begin{equation}
\begin{split}
&\int_0^{T^*}\intO{ \kappa(\thet)\nabla \mathcal{T}_{M, \delta}(\thet) \cdot \nabla \mathcal{T}_k \left(\mathcal{T}_{M, \delta}(\thet) - \thetah\right)}\dt\\
&\ge-\intO{\mathcal{G}_k\left( \mathcal{T}_{M, \delta}(\thet)(T^*)- \thetah\right)-\mathcal{G}_k\left( \mathcal{T}_{M, \delta}(\thet_0)- \thetah\right)}\\
& + \int_0^{T^*}\intO{\left(\mathcal{T}^{\,'}_{M, \delta}(\thet) \, \S:\D{\vv}-\nabla [\mathcal{T}_{M, \delta}(\thet)]\cdot \vv   \right) \mathcal{T}_k \left( \mathcal{T}_{M, \delta}(\thet) - \thetah\right)}\dt\\
&-k\left(\int_{\{\thet_0>M/2\}}\!\!\!\!\! \thet_0\dx + \int_{\{\thet>M/2\}}\!\!\!\!\!\S:\D{\vv} \dx\dt\right).
\end{split}\label{ep}
\end{equation}
For $M> k+1 + \|\thetah\|_{\infty}$, we have $\mathcal{T}_k (  \mathcal{T}_{M, \delta}(\thet)  - \thetah) = \mathcal{T}_k ( \thet  - \thetah)$, thus one can let $M\to \infty$ in \eqref{ep}. Using integration by parts and the fact that $\diver \vv=0$ similarly as in \eqref{thetab-n} one obtains % $ \nabla \mathcal{T}_{M, \delta}(\thet) \cdot \nabla \mathcal{T}_k \left(  \mathcal{T}_{M, \delta}(\thet)  - \thetah\right) = \nabla \thet \cdot \nabla  \mathcal{T}_k \left( \thet  - \thetah\right)$ and $ \mathcal{T}^{\,'}_{M, \delta}(\thet) \,  \mathcal{T}_k \left(  \mathcal{T}_{M, \delta}(\thet)  - \thetah\right) =  \mathcal{T}_k \left( \thet  - \thetah\right)$.
\begin{equation*}
\begin{split}
&\int_0^{T^*}\intO{ \kappa(\thet)\nabla \thet \cdot \nabla \mathcal{T}_k (\thet - \thetah)}\dt\\
&\ge-\intO{\mathcal{G}_k (\thet(T^*)- \thetah)-\mathcal{G}_k(\thet_0- \thetah)} + \int_0^{T^*}\intO{( \S:\D{\vv}-\nabla \thetah\cdot \vv ) \mathcal{T}_k ( \thet - \thetah)}\dt.
\end{split}%\label{ep}
\end{equation*}

%Moreover, we know from \eqref{T2} that
%$$
%\int_0^{T^*}\intO{ \fint_t^{t +\ep} k(\thet^n) \mathcal{T}^{\,''}_{M, \delta}(\thet^n) |\nabla\thet^n|^2\dtau\ \mathcal{T}_k \left(\mathcal{T}_{M, \delta}(\thet^n)- \thetah\right)}\dt\to0\mbox{ as $M\to+\infty$.}
%$$
%Taking $M\to+ \infty$ in \eqref{ep} we obtain
%\begin{equation}\begin{split}
%&\intO{\mathcal{G}_k\left(\thet(T^*)- \thetah \right)} -\intO{\mathcal{G}_k\left(  \thet_0  - \thetah\right)}
%+ \int_0^{T^*}\intO{ \nabla \thetah \cdot \vv\ \mathcal{T}_k \left( \thet  - \thetah\right)}\dt \\&
% - \int_0^{T^*}\intO{   \mathcal{T}_k \left( \thet  - \thetah\right) \ \S:\D{\vv}}\dt
%=-  \int_0^{T^*}\intO{ k(\thet) \nabla \thet \cdot \nabla \mathcal{T}_k \left(  \thet  - \thetah\right)}\dt.
%\end{split}\end{equation}
Comparing the result with \eqref{MB}, we see that
\begin{equation}\label{GradTk-L2}
\begin{split}
  \limsup_{n\to +\infty} %&
  \int_0^{T^*}\intO{\kappa(\thet^n)\nabla \thet^n \cdot \nabla \mathcal{T}_k(\thet^n -\thetah)}\dt
  %\\&
  \le  \int_0^{T^*}\intO{ \kappa(\thet) \nabla \thet \cdot \nabla \mathcal{T}_k (\thet  - \thetah )}\dt,
\end{split}
\end{equation}
which is the corner stone for the strong convergence. Indeed, it follows from \eqref{thet-weak}--\eqref{thet-strong} and from \eqref{GradTk-L2} that (using also the Vitali convergence theorem and \eqref{k})
\begin{equation}\label{goal}
\begin{split}
&\limsup_{n\to +\infty} \int_0^{T^*}\intO{\kappa(\thet^n)|\nabla \mathcal{T}_k(\thet^n -\thetah)|^2}\dt \le  \int_0^{T^*}\intO{ \kappa(\thet) |\nabla \mathcal{T}_k (\thet  - \thetah )|^2}\dt.
\end{split}
\end{equation}
Since
\begin{align}\label{weakGradTk}
\sqrt{\kappa(\thet^n)}\nabla \mathcal{T}_k(\thet^n -\thetah) &\rightharpoonup \sqrt{\kappa(\thet)}\nabla \mathcal{T}_k(\thet -\thetah) &&\mbox{ weakly in } L^2(Q),\\
\intertext{the inequality \eqref{goal}, weak lower semicontinuity of $L^2$ norm and uniform convexity of $L^2$ imply} %therefore \eqref{GradTk-L2} becomes
%\begin{equation}\begin{split}
%\lim_{n\to +\infty} \int_0^{T^*}\intO{&k(\thet^n)\nabla (\thet^n - \thetah) \cdot \nabla [\mathcal{T}_k(\thet^n -\thetah)]}\dt
%\\ &= \int_0^{T^*}\intO{ k(\thet) \nabla (\thet-\thetah) \cdot \nabla [\mathcal{T}_k \left(  \thet  - \thetah\right)]}\dt,
%\end{split}\end{equation}
%The fact that $\mathcal{T}^{\,'}_k(s)=0$ for $s>k$ yields that
% \begin{equation}\label{norm}
%\lim_{n\to +\infty} \int_0^{T^*}\intO{k(\thet^n) |\nabla [\mathcal{T}_k(\thet^n -\thetah)]|^2}\dt= \int_0^{T^*}\intO{ k(\thet) |\nabla [\mathcal{T}_k \left(  \thet  - \thetah\right)]|^2}\dt.
%\end{equation}
%Weak convergence \eqref{weakGradTk} with $\alpha=1/2$ combined with \eqref{norm} and uniform convexity of $L^2((0,T^*)\times\Omega)$ gives
\label{strong-Tk}
\sqrt{\kappa(\thet^n)}\nabla \mathcal{T}_k(\thet^n -\thetah) &\to \sqrt{\kappa(\thet)}\nabla \mathcal{T}_k(\thet -\thetah) &&\mbox{ strongly in } L^2((0, T^*)\times\Omega).
\end{align}

Finally, we want to show that
\begin{align}\label{grad-eta}
\sqrt{\kappa(\thet^n)}\frac{\nabla \thet^n}{\thet^n} \to \sqrt{\kappa(\thet)} \frac{\nabla \thet}{\thet} &&\mbox{ strongly in } L^2((0, T^*)\times\Omega).
\end{align}
We have
\begin{align}\label{weak-nablathet-thet}
\sqrt{\kappa(\thet^n)}\frac{\nabla \thet^n}{\thet^n} \rightharpoonup \sqrt{\kappa(\thet)}\frac{\nabla \thet}{\thet} &&\mbox{ weakly in } L^2((0, T^*)\times\Omega).
\end{align}
because of the uniform boundedness in $L^2((0, T^*)\times\Omega).$
%similarly as in \eqref{weakGradTk} with help of Vitali theorem, \eqref{conv-point}, \eqref{ae-nabla} and \eqref{uniform-alpha}.
%
To show convergence of norms, we write
\begin{equation}\label{nablathet-thet}
\begin{split}
& \int_0^{T^*}\intO{ \kappa(\thet^n)\frac{|\nabla \thet^n|^2}{(\thet^n)^2}}\dt \\
 &= \int_0^{T^*}\intO{\kappa(\thet^n)\frac{\nabla \thet^n}{\thet^n}\cdot\frac{\nabla (\thet^n - \thetah)}{\thet^n}}\dt +  \int_0^{T^*}\intO{\kappa(\thet^n)\frac{\nabla \thet^n}{\thet^n}\cdot\frac{\nabla \thetah}{\thet^n}}\dt\\
&= \int_0^{T^*}\intO{\kappa(\thet^n)\frac{\nabla \thet^n}{\thet^n}\cdot\frac{\nabla \mathcal{T}_k(\thet^n - \thetah)}{\thet^n}}\dt + \int_{\{|\thet^n-\thetah|>k\}}\!\!\!\!\! \kappa(\thet^n) \frac{\nabla \thet^n}{\thet^n}\cdot\frac{\nabla (\thet^n - \thetah)}{\thet^n} \dx\dt \\&+
 \int_0^{T^*}\intO{\kappa(\thet^n)\frac{\nabla \thet^n}{\thet^n}\cdot\frac{\nabla \thetah}{\thet^n}}\dt.
\end{split}
\end{equation}
Next, we let $n\to \infty$ on the right-hand side of \eqref{nablathet-thet}. First, by the weak convergence \eqref{weak-nablathet-thet}, the strong convergence results \eqref{thet-strong} and \eqref{strong-Tk}, the Vitali theorem and the minimum principle $\thet^n\ge \mu$
\begin{equation*}
\lim_{n\to + \infty}\int_0^{T^*}\!\!\!\intO{\kappa(\thet^n)\frac{\nabla \thet^n}{\thet^n}\cdot\frac{\nabla \mathcal{T}_k(\thet^n - \thetah)}{\thet^n}}\dt = \int_0^{T^*}\!\!\!\intO{\kappa(\thet)\frac{\nabla \thet}{\thet}\cdot\frac{\nabla \mathcal{T}_k(\thet - \thetah)}{\thet}}\dt
\end{equation*}
and
\begin{equation} \lim_{n\to+\infty}  \int_0^{T^*}\intO{\kappa(\thet^n)\frac{\nabla \thet^n}{\thet^n}\cdot\frac{\nabla \thetah}{\thet^n}}\dt =  \int_0^{T^*}\intO{\kappa(\thet)\frac{\nabla \thet}{\thet}\cdot\frac{\nabla \thetah}{\thet}}\dt.
\end{equation}
In addition, assuming that $k>\|\thetah\|_{\infty}$, we see that $\thet^n\geq k$ on the set $\{|\thet^n-\thetah|>k\}$. Consequently, with the help of the Young inequality,  we deduce that for any $\lambda\in(0,1)$
\begin{equation}\label{limit-k}
\begin{split}
\left|\int_{\{|\thet^n-\thetah|>k\}}\!\!\!\!\! \kappa(\thet^n)\frac{\nabla \thet^n}{\thet^n}\cdot\frac{\nabla (\thet^n - \thetah)}{\thet^n} \dx\dt \right|& \leq  2\overline{\kappa}\intTO{\frac{|\nabla \thet^n|^2+ |\nabla\thetah|^2}{(\thet^n)^{1+\lambda} k^{1-\lambda} }}\leq \frac{C(\lambda)}{k^{1-\lambda}},
\end{split}
\end{equation}
where $C$ is independent of $n$ and $k$ thanks to the uniform estimate in \eqref{eq:re3} and the minimum principle $\thet^n\ge \mu$.
Equality \eqref{nablathet-thet} and inequality \eqref{limit-k} remain valid also if we replace $\thet^n$ with $\thet$.
%In the exactly same manner, one can replace $\thet^n$ by  $\thet$ in~\eqref{nablathet-thet}--\eqref{limit-k}.
Consequently, it follows%choosing for a fixed $\varepsilon>0$ suitable $k>0$ and then $n_0\in\N$ we get for $n>n_0$ that
\begin{equation}\label{conv-norms}
\limsup_{n\to \infty}\int_0^{T^*}\intO{\kappa(\thet^n)\frac{|\nabla \thet^n|^2}{(\thet^n)^2}}\dt \le \int_0^{T^*}\intO{\kappa(\thet)\frac{|\nabla \thet|^2}{\thet^2}}\dt +\frac{C(\lambda)}{k^{1-\lambda}}.
\end{equation}
Since $k>\|\thetah\|_{\infty}$ is arbitrary, \eqref{conv-norms} gives convergence of norms which combined with \eqref{weak-nablathet-thet} implies the strong convergence \eqref{grad-eta}. Finally,  let us recall that $T^*$ can be chosen arbitrarily close to $T$. In addition, we can a~priori construct the solution on the time interval $(0,2T)$. Consequently, $T^*$ can be chosen bigger than $T$ and so we obtain \eqref{grad-eta} with $T^*$ replaced by $T$.

\subsubsection{Limit in the entropy equality}

Let us analyze the limit of \eqref{entropy-weak} as $n\to+\infty$ for each term.
The strong convergence of $\thet^n$ and $\nabla\thet^n$, see \eqref{thet-strong} and \eqref{NS1}, implies that $\eta^n\to\eta$ and $\nabla\eta^n\to\nabla\eta$ a.e. in $Q$. Thanks to the definition of $\eta^n$, and the a~priori bound \eqref{eq:re3}, we obtain a uniform bound for $\nm{\eta^n}_{L^r(Q)}$ for any $r\in(1,+\infty)$ and for $\nm{\nabla\eta^n}_{L^2(Q)}$. These facts, the Vitali convergence theorem and \eqref{thet-zero} allow us, up to subsequence, easily pass to the limit on the left hand side of~\eqref{entropy-weak}.

To pass to the limit also in terms on the right hand side, we use \eqref{Lp}, \eqref{Lp-S}, \eqref{thet-strong}, \eqref{Mintye} and \eqref{grad-eta} combined with the Lebesgue dominated convergence theorem and together with the minimum principle for $\thet^n$.  %Since $\S^n:\D{\vv^n}\to \S:\D{\vv}$ in $L^1((0,T)\times\Omega)$, we can assume that there exists an integrable majorant for a subsequence, denoted same, of $\{\S^n:\D{\vv^n}\}$. This together with \eqref{conv-point} allows to pass to the limit in the first term on the right side of \eqref{entropy-weak} by Dominated Convergence Theorem. The limit passage in the last term is clear due to \eqref{conv-norms}.
Consequently, $\eta$ satisfies entropy equation~\eqref{entropy-limit}.

\subsection{Continuity of \texorpdfstring{$\thet$}{theta} in time}

Finally, we focus on the attainment of initial conditions and continuity with respect to time variable. Concerning the velocity field, we can recall \eqref{Lp} and by standard parabolic interpolation, we observe $\vv \in C([0,T];  L^2_{0,\diver})$. The fact that $\vv(0)=\vv_0$ is then proven analogously as for Navier--Stokes equations, see e.g. \cite{MaBook,FrMaRu10}.

Now, we focus on the temperature. For the attainment of the initial condition, one can follow \cite{Consiglieri,BuMaRa09,FrMaRu10}. Thus we present here only the proof of continuity of $\thet$ with respect to time into the $L^1$ topology. The key idea is the following. We investigate a function (here $M>\max\{\|\thetah\|_{L^{\infty}(\Omega)},2,2\mu\}$)
$$
g(\thet-\thetah):= \sign (\thet - \thetah) \sqrt{\mathcal{G}_M(\thet - \thetah)}
$$
and show first that for all $\psi \in L^2(\Omega)$ the mapping $t\mapsto \intO{g(\thet(t,x)-\thetah(t,x))\psi(x)}$ is continuous on $[0,T]$. Second, we show that the mapping $t\mapsto \intO{g^2(\thet(t,x)-\thetah(t,x))}$ is continuous on $[0,T]$. As a direct consequence of these two properties, we obtain that
\begin{equation}
g(\thet-\thetah)\in C([0,T]; L^2(\Omega)).
\label{C2}
\end{equation}
Since,
$$
g'(\thet - \thetah)=\frac{|T_M(\thet - \thetah)|}{2\sqrt{\mathcal{G}_M(\thet - \thetah)}}= \left\{
\begin{aligned}
&\frac{1}{\sqrt{2}},&&\textrm{if }\thet\le M+\thetah\\
&\frac{M}{2\sqrt{M(\thet - \thetah)-\frac{M^2}{2}}}&&\textrm{if }\thet\ge M+\thetah
\end{aligned}
\right\} \ge \frac{\sqrt\mu}{\sqrt{2\thet}},
$$
we have a trivial estimate
$$
|\thet_1-\thet_2|\le \sqrt{2}(\sqrt{\thet_1}+\sqrt{\thet_2})|g(\thet_1 - \thetah)-g(\thet_2 - \thetah)|.
$$
Consequently, using the H\"{o}lder inequality and \eqref{theinf1}, we observe
$$
\begin{aligned}
\|\thet(t_1)-\thet(t_2)\|_{L^1(\Omega)} &\le \sqrt{2}\|\sqrt{\thet(t_1)}+\sqrt{\thet(t_2)}\|_{L^2(\Omega)} \|g(\thet(t_1) - \thetah)-g(\thet(t_2) - \thetah)\|_{L^2(\Omega)}\\
&\le C \|g(\thet(t_1) - \thetah)-g(\thet(t_2) - \thetah)\|_{L^2(\Omega)}.
\end{aligned}
$$
Thus, we see that \eqref{C2} implies $\thet \in C([0,T];L^1(\Omega))$.
% In addition setting $\varphi:= (1-t/\delta)\chi_{(0,\delta)}\psi(x)$ with $\delta\in (0,T)$ and arbitrary $\psi \in C_0^1(\Omega)$, we get
% \begin{equation*}
% \begin{aligned}%\label{T2}
% \fint_0^{\delta}\intO{\thet(t) \psi}\dt = \int_0^{\delta}(1-t/\delta)\intO{\thet\vv\cdot\nabla\psi- \kappa(\thet)\nabla\thet\cdot\nabla\psi}\dt \\
%  + \int_0^{\delta}(1-t/\delta)\intO{\S:\D\vv\, \psi}\dt + \intO{\thet_0\psi} %\mbox{ for all } \varphi\in C^\infty_0((-\infty, T); C_0^\infty(\Omega))
% \end{aligned}
% \end{equation*}
% and then letting $\delta \to 0_+$ and using also $\thet \in C([0,T];L^1(\Omega))$, we see
% \begin{equation*}
% \begin{aligned}%\label{T2}
% \intO{\thet(0) \psi}=\lim_{\delta\to 0_+}\fint_0^{\delta}\intO{\thet(t) \psi}\dt = \intO{\thet_0\psi}. %\mbox{ for all } \varphi\in C^\infty_0((-\infty, T); C_0^\infty(\Omega))
% \end{aligned}
% \end{equation*}
% Since $\psi$ is arbitrary, we have $\thet(0)=\thet_0$.

It remains to show \eqref{C2}. First, we set $\psi:=g'(\thet^n-\thetah)\chi_{[0,\tau]}  \varphi$ with $\varphi \in C^1_0(\Omega)$ arbitrary in \eqref{ode22}. Using integration by parts, we have
\begin{equation*}
\begin{aligned}
%\label{ode22n}
&\intO{(g(\thet^n(\tau)-\thetah)-g(\thet^n(0)-\thetah))\varphi}= -\int_0^{\tau}\intO{g'(\thet^n-\thetah) \varphi \vv^{n} \nabla \thetah}\dt\\
&  + \int_0^{\tau}\intO{\left(\vv^{n} (\thet^{n}-\thetah)-\kappa(\thet^n)\nabla \thet^n\right) \cdot (\varphi g''(\thet^n-\thetah) \nabla (\thet^n-\thetah) +g'(\thet^n-\thetah) \nabla \varphi)}\dt\\
& + \int_0^{\tau}\intO{\S^{n}: \D\vv^{n} g'(\thet^n-\thetah)  \varphi} \dt.
\end{aligned}
\end{equation*}
Next, we let $n\to \infty$ in the above identity. Using the uniform estimates \eqref{apr:vn}--\eqref{eq:re3}, the properties of $g$, the convergence results \eqref{Lp}, \eqref{thet-strong}, \eqref{ae-nabla}, \eqref{weak-nablathet-thet} and the Vitali convergence theorem, we can easily identify the limits in the first two terms on the right hand side. For the last term on the right hand side, we also use \eqref{Sn}. In addition, for almost all $\tau \in (0,T)$ we can also identify the limit on the left hand side
\begin{equation}
\begin{aligned}
\label{ode22n}
&\intO{(g(\thet(\tau)-\thetah)-g(\thet_0-\thetah))\varphi}= -\int_0^{\tau}\intO{g'(\thet-\thetah) \varphi \vv \nabla \thetah}\dt\\
&  + \int_0^{\tau}\intO{\left(\vv (\thet-\thetah)-\kappa(\thet)\nabla \thet\right) \cdot (\varphi g''(\thet-\thetah) \nabla (\thet-\thetah) +g'(\thet-\thetah) \nabla \varphi)}\dt\\
& + \int_0^{\tau}\intO{\S: \D\vv g'(\thet-\thetah)  \varphi} \dt.
\end{aligned}
\end{equation}
Since the right hand side is a continuous function of $\tau \in [0,T]$, we can redefine $g(\thet-\thetah)$ on zero subset of $[0,T]$ to get
\begin{equation}\label{weak-cont1}
\left(\intO{g(\thet)\varphi}\right) \in C([0,T]) \qquad \textrm{ for all }\varphi \in C_0^1(\Omega).
\end{equation}
Since $g\in L^{\infty}(0,T; L^2(\Omega))$ and $C_0^1(\Omega)$ is dense in $L^2(\Omega)$, \eqref{weak-cont1} implies that
\begin{equation}\label{WC2}
\left(\intO{g(\thet)\varphi}\right) \in C([0,T]) \qquad \textrm{ for all }\varphi \in L^2(\Omega).
\end{equation}

Finally, we pass in \eqref{MBa} to the limit as $n\to+\infty$ similarly as in \eqref{MB} using also \eqref{weak-nablathet-thet}
\begin{equation}\label{MBaa}
\begin{split}
&\intO{\mathcal{G}_M(\thet(\tau) -\thetah)-\mathcal{G}_M(\thet_0 -\thetah)}=\int_0^{\tau}\intO{-\kappa(\thet)\nabla \thet \cdot \nabla [\mathcal{T}_M(\thet -\thetah)]} \dt\\
& + \int_0^{\tau}\intO{-\nabla \thetah \cdot \vv \, \mathcal{T}_M(\thet -\thetah)+\mathcal{T}_M(\thet -\thetah)\, \S:\D{\vv} } \dt=:\int_0^{\tau}\intO{h}\dt,
\end{split}
\end{equation}
where $h\in L^{1}(Q)$. Hence, we see that (it can be continuously extended)
$$
\intO{\mathcal{G}_M(\thet(\tau) -\thetah)} \in C([0,T]).
$$
Since $(g(\theta-\thetah))^2= \mathcal{G}_M(\thet -\thetah)$, the above relation combined with \eqref{WC2} implies \eqref{C2}. The proof is\footnote{In fact in above procedure we somehow extended $g$ and $\mathcal{G}_M$ also to a possible non-Lebesgue points. This can be done more carefully. Namely, one can consider $\fint_t^{t+\delta} g(\tau)\dtau$ and $\fint_{t}^{t+\delta}\mathcal{G}_M(\thet(\tau)-\thetah)$. These are surely continuous with respect to $t\in [0,T]$. Then thanks to \eqref{WC2} and \eqref{MBaa}, we see that
$$
\fint_t^{t+\delta} g(\tau)\dtau \to g \quad \textrm{strongly in } C([0,T]; L^2(\Omega)).
$$
} complete. \qed

\section*{Acknowledgment}
{A.~Abbatiello has been supported by the ERC-STG Grant n. 759229 HiCoS ``Higher Co-dimension Singularities: Minimal Surfaces and the Thin Obstacle Problem" and is member of the Italian National Group for the Mathematical Physics (GNFM) of the Italian National Institute of the High Mathematics (INdAM).   M. Bul\'{\i}\v{c}ek and P. Kaplick\'{y} acknowledge the support of the project  No. 20-11027X financed by Czech Science Foundation (GA\v{C}R). M. Bul\'{\i}\v{c}ek is member of the Jind\v{r}ich Ne\v{c}as Center for Mathematical Modelling.}
\bibliographystyle{amsplain}
\providecommand{\bysame}{\leavevmode\hbox to3em{\hrulefill}\thinspace}
\providecommand{\MR}{\relax\ifhmode\unskip\space\fi MR }
\providecommand{\MRhref}[2]{%
  \href{http://www.ams.org/mathscinet-getitem?mr=#1}{#2}
}
\providecommand{\href}[2]{#2}

\end{document}